%% file: NCAndHnCenters.tex
\documentclass[12pt]{article}      
\usepackage{makeidx}
\usepackage[left=0.5in, right=0.5in, top=0.5in, bottom=0.5in, textwidth=2in]{geometry} 
\usepackage{amsmath, amsfonts, epsfig,amssymb}
\usepackage[all]{xy}

\input macros-II.tex

\input diagrams-quotient.tex

\newtheorem{claim}[thm]{Claim}

\newcommand{\M}{M\"{o}bius band}
\newcommand{\AT}{\ensuremath{\mathcal{A}/[\mathcal{A},\mathcal{A}]_t}}
\newcommand{\A}{\ensuremath{\mathcal{A}}}
\newcommand{\AC}{\ensuremath{\mathcal{A}/[\mathcal{A},\mathcal{A}]}}
\renewcommand{\NC}{\ensuremath{\mathcal{NC}_n}}
\newcommand{\Sn}{\ensuremath{\mathcal{S}_n}}
\newcommand{\Hn}{\ensuremath{\mathcal{H}_n}}
\renewcommand{\C}{\ensuremath{\mathbb{C}}}
\newcommand{\real}{\ensuremath{\mathbb{R}}}

\title{The Center of the Nilcoxeter and $0$-Hecke Algebras}  \author{Joelle Brichard}     \date{\today }      

\begin{document}             
\maketitle   
\begin{abstract}
We compute the dimension of the center of the 0-Hecke algebra \Hn{} and of the Nilcoxeter algebra \NC{} using a calculus of diagrams on the \M.  In the case of the Nilcoxeter algebra, this calculus is shown to produce a basis for $Z(\NC)$ and the table of multiplication in this basis is shown to be trivial.  We conjecture that a basis for $Z(\Hn)$ can also be obtained in a specific  way from this topological calculus.
\end{abstract}             

\section{Introduction}

	\subsection {Generic Algebras}
Given a Coxeter system $(W,S)$ and an associative commutative algebra $\A$ with unity, we start by defining the generic algebra $\cal{E}  _{\A} (a_s,b_s)$. This approach is taken from \cite{Humphrey}
and depends on parameters $\{a_s, b_s \in \A \}_{s\in S}$  subject to the condition that $a_s=a_t$ and $b_s=b_t$ whenever $s$ and $t$ are conjugates in $W$.   $\cal{E} _{\A} (a_s,b_s)$ is generated by elements $T_w$ for $w\in W$ with relations
 
\begin{enumerate}
\item
$T_sT_w=T_{sw}\;\;\;$ if $l(sw)>l(w)$
\item
$T_sT_w=a_sT_w+b_sT_{sw}\;\;\;$ if $l(sw)<l(w).$
\end{enumerate}

The two algebras of concern in this paper are both examples of this type of generic algebras with Coxeter group \Sn.  Other important examples include the group algebra $\A[W]$ which is obtained by setting all $a_s=0$ and $b_s=1$.  Hecke algebras are also examples of this construction for $\A=\mathbb{Z} [q,q^{-1}]$ and with $a_s=q-1$ and $b_s=q$ for all $s\in S$.  For more on Coxeter groups and Hecke algebras see \cite{KL}.

There is another set of relations equivalent to the ones above:
\begin{enumerate}
\item
$T_sT_w=T_{sw}\;\;\;$ if $l(sw)<l(w),$
\item
$T_s^2=a_sT_s+b_sT_e$
\end{enumerate}
with $T_e$ the identity element.  $\cal{E}$ is in fact generated as an algebra by the $T_s,\, s\in S$ and $1=T_e$.

Any generic algebra with coxeter group \Sn, is generated by the following three relations:
\begin{enumerate}
\item
$T_iT_{i+1}T_i=T_{i+1}T_iT_{i+1}$
\item
$T_iT_j=T_jT_i$ for $|i-j|>1$
\item
$T_i^2=a_iT_i+b_iT_e$
\end{enumerate}

		\subsubsection{The Nilcoxeter Algebra} \label{NC}
			
The Nilcoxeter algebra, which we denote \NC, is an example of a generic algebra associated to  \Sn{}  where we set $a_s=b_s=0$ for all $s\in S$.   
This algebra first appeared in \cite{BGG} in relation to the cohomology of flag varieties.  It was also studied by Lascoux and Sch\"{u}tzenberger \cite{LS}, Macdonald \cite{Mc}, Fomin and Stanley \cite{FS}
and others.  Khovanov has shown in \cite{K:Nilcoxeter} that it categorifies the polynomial representation of the Weyl algebra and recent work by Khovanov and Lauda showed its relevance to the categorification of quantum groups \cite{KhL}.  
In terms of the generators $T_i,$ where $T_i=T_{s_i}$ for $s_i\in\Sn$ the transposition $(i\;\;i+1),$ we have the following defining relations:
\begin{enumerate}
\item
$T_iT_{i+1}T_i=T_{i+1}T_iT_{i+1}$
\item
$T_iT_j=T_jT_i$ for $|i-j|>1$
\item
$T_i^2=0$
\end{enumerate}

We therefore also have
$$T_iT_w=\left\{ \begin{array}{ll}T_{s_iw}&if \,\,l(s_iw)>l(w) \\0&if\,\, l(s_iw)<l(w)\end{array}\right. .$$

		\subsubsection{The 0-Hecke Algebra} \label{Hn}

The 0-Hecke algebra, denoted \Hn, is the Hecke algebra associated to  \Sn{} with $q=0$.  Equivalently, it is the generic algebra where we have set $a_s=1$ and $b_s=0$ for all $s\in S$.  That is, the third relation on the $T_i$ now becomes $$T_i^2=T_i.$$  This implies that 
$$T_iT_w=\left\{ \begin{array}{ll}T_{s_iw}&if \,\,l(s_iw)>l(w) \\
T_w&if\,\, l(s_iw)<l(w)\end{array}\right.$$
Norton studied the Hecke algebra of a Coxeter group at $q=0$ for all types in \cite{Norton}; she classified the irreducible modules, described the decomposition of the algebra into left ideals as well as the Cartan invariants.  In \cite{Carter}, Carter gave decomposition numbers for \Hn{} in the type $A$ case.  Krob and Thibon gave a representation-theoretic interpretation of non-commutative algebras \cite{KT}
 and works by Duchamp, Hivert and Thibon \cite{DHT} showed that the representations of \Hn{} in type $A$ categorify the ring of quasi-symmetric functions.

	\subsection{The Centers of Frobenius Algebras}

We now derive some general facts about the centers of Frobenius algebras and then introduce a trace map which makes both \NC{} and \Hn{} into Frobenius algebra.  The facts here derived are crucial to our approach to describing the centers of \NC{} and \Hn.

		\subsubsection{Symmetric Frobenius Algebras}

To investigate the centers of Frobenius algebras, it is sometimes convenient to use a duality relation between the center $Z(\A)$ and a quotient of  $\A$ by certain subgroups.  In the case of symmetric Frobenius algebras, this subgroup is simply the commutator subgroup.
Let $\A$ be a Frobenius algebra over a field $k$, $\epsilon : \A \rightarrow k$ its trace map and assume that $\A$ is a symmetric algebra, meaning that $\epsilon(ab)=\epsilon(ba)$ for all $ a,b$.  Denote by $\A^*=$Hom$_k(\A,k)$ the dual of $\A$. Then $\epsilon$ extends to an isomorphism of $\A$-bimodules: 
\begin{eqnarray*}
\phi:\A&\rightarrow& \A^*\\
1&\mapsto&  \epsilon\\
a&\mapsto& \epsilon( \,*a).
\end{eqnarray*}
Let  $[\A,\A]$ denote the commutator subspace of $\A$: the space spanned by elements of the form $ab-ba$ for $ a,b \in \A.$\\
\\
\begin{claim} If $\A$ is a symmetric Frobenius algebra, then $$(\A/[\A,\A])^*\simeq Z(\A).$$
\end{claim}
\begin{proof}
Let 
\begin{eqnarray*}
\Psi: Z(\A)&\rightarrow& (\A/[\A,\A])^*\\
z&\mapsto &z\epsilon :\A\rightarrow k; z\epsilon(a)=\epsilon(za)=\epsilon(az).
\end{eqnarray*}
$\Psi $ is well-defined: $$z\epsilon ([a,b])=z\epsilon (ab-ba)=\epsilon (zab-zba)=\epsilon (bza)-\epsilon (zba)=\epsilon(zba)-\epsilon(zba)=0$$ since $z$ is a central element and $\A$ is
symmetric.  \\
$\Psi$ is injective: since the map $\A \rightarrow \A^*$ is an isomorphism and the map $Z(\A)\rightarrow \A$ is an inclusion, the composition $ Z(\A)\rightarrow \A\rightarrow \A^*$ is injective.
Moreover, we have seen that the map $Z(\A)\rightarrow \A^*$ factors through $Z(\A)\rightarrow (\A/[\A,\A])^*$, which means that $\Psi : Z(\A)\rightarrow (\A/[\A,\A])^*$ is
injective as well.\\
$\Psi $ is surjective: let $\alpha \in (\A/[\A,\A])^*$.  Because $$\pi : \A\rightarrow \A/[\A,\A]$$ is surjective, its dual $$\pi ^* (\A/[\A,\A])^*\rightarrow \A^*$$ must also be
injective.  We know from the isomorphism $\phi:\A\rightarrow \A^*$ that $\pi^*(\alpha)=z\epsilon$ for some $z\in \A$.  We must check that $z\in Z(\A)$.  Now,
since $\alpha \in (\A/[\A,\A])^*$, $z\epsilon([a,b])=0 $ for all $ a,b \in \A.$  Hence $$\epsilon(z[a,b])=\epsilon(zab)-\epsilon(zba)=\epsilon(zab)-\epsilon(azb)=\epsilon([z,a]b)=0.$$  Since the trace $\epsilon$
is non-degenerate, this cannot be true for all $b \in \A$ unless $[z,a]=0$.  This is true for all $a\in \A$ as well, so we conclude that $z\in Z(\A)$, with 
$\Psi(z)=\alpha$.

\end{proof}

		\subsubsection{Non-Symmetric Frobenius Algebras}
We now look at the case when $\A$ is not a symmetric Frobenius algebra.  Since $$\epsilon:\A\to k$$ is non-degenerate, one can define an involution $f:\A\to \A$ such that $\epsilon(ab)=\epsilon(bf(a)),\,f^2(a)=a,\,f(ab)=f(a)f(b)$.   Using this involution, one can state a duality very similar to the one we have derived in the case of the symmetric Frobenius algebras.

Given the involution $f$ on $\A$, let us denote the "twisted" commutator of $a$ and $b$ by $[a,b]_t=ab-bf(a)$.  Similarly, let $[\A,\A]_t$ be the subalgebra 
of $\A$ spanned by elements $ab-bf(a)$.  We also define the twisted center of \A{} as follows: 
$$TZ(\A)=\{b\in \A|ba-f(a)b=0 \text{ for all } b \in \A\}.$$ 
Note that if $b\in TZ(\A)$ then $f(b)\in TZ(\A)$ as well because 
$$f(f(b)a)=bf(a)=at=f(f(a)f(b)),$$ which means that $f(b)a=f(a)f(b)$ since the involution $f$ is bijective. Moreover, if $b\in TZ(\A)$, we have 
$$f(b)f(a)=f(ba)=f(f(a)b)=af(b),$$ such that for a general $t\in TZ(\A)$, $at-tf(a)=0$.
\\
\begin{claim} \label{TwistedDuality}
With the above notation,\\
 \begin{eqnarray}
 Z(\A)&\backsimeq & (\A/[\A,\A]_t)^*\\ 
 TZ(\A) &\backsimeq& (\A/[\A,\A])^*.
\end{eqnarray}
\end{claim}
\begin{proof} (1) The proof closely resembles the one for $\A$ symmetric.  Let 
\begin{eqnarray*} 
\Psi: Z(\A)&\rightarrow& (\A/[\A,\A]_t)^*\\
z&\mapsto& z\epsilon :\A\rightarrow k; z\epsilon(a)=\epsilon(za).
\end{eqnarray*}
$\Psi $ is well-defined: 
$$z\epsilon ([a,b]_t)=z\epsilon (ab-bf(a))=\epsilon (zab-zbf(a))=\epsilon (a(zb))-\epsilon ((zb)f(a))=0$$ 
since $z$ is central and $\epsilon(ab)=\epsilon(bf(a)).$\\
$\Psi $ is injective: The proof is exactly as for $\A$ symmetric.\\
$\Psi$ is surjective: Let $\alpha\in (\A/[\A,\A]_t)^*$ and $\pi^*:(\A/[\A,\A]_t)^*\rightarrow \A^*$.  Then $\pi^*(\alpha)=z\epsilon$ for some $z\in \A$ and $z\epsilon([a,b]_t)=0 $ for all $
a,b\in \A.$  We must verify that $z$ is in the center $Z(\A)$.\\
$$0=z\epsilon([a,b]_t)=\epsilon(zab)-\epsilon(zbf(a))=\epsilon(zab)-\epsilon(azb)=\epsilon([z,a]b).$$
 Since this is true for all $b\in \A$ and since $\epsilon$ is non-degenerate, then $[z,a]$ must vanish for
all $a\in \A$, which means that $z$ is central.\\
(2) Let 
\begin{eqnarray*}
\Psi: TZ(\A)&\rightarrow& (\A/[\A,\A])^*\\
t\mapsto t\epsilon :\A&\rightarrow& k; t\epsilon(a)=\epsilon(ta).
\end{eqnarray*}
$\Psi $ is well-defined:
$$c\epsilon([a,b])=\epsilon(cab)-\epsilon(cba)=\epsilon(cab)-\epsilon(baf(c))=\epsilon(cab)-\epsilon(bf(c)f(a))=\epsilon(cab)-\epsilon(bf(ca))=0.$$
$\Psi$ is injective: The proof is as before.\\
$\Psi$ is surjective:  We have as before $c\epsilon=\pi^*(\alpha)$ and we need to check that $c\in TZ(\A).$ $$0=c\epsilon([a,b]_t)=\epsilon(cab)-\epsilon(cba)
=\epsilon(cab)-\epsilon(f(a)cb)=\epsilon([t,a]_tb)$$ so that $[c,a]_t=0$ and $c\in TZ(\A)$.
\end{proof}

Note that $\A/[\A,\A]$ and \AT are no longer algebras, but vector spaces, both for $\A=\NC$ and $\A=\Hn$.

		\subsubsection{The Nilcoxeter and 0-Hecke Algebras as Non-Symmetric Frobenius Algebras} 

We can endow the $0$-Hecke and Nilcoxeter algebras with a trace map, making them Frobenius algebras. Some of our constructions and arguments will be similar for both \NC{} and \Hn. In these cases we denote either of them by $\A_n$, where $\A_n$ is generated by $T_1, T_2, \dots ,T_{n-1}$, if we want to emphasize the dimension of our algebra, or simply by \A. Then $\epsilon:\A\rightarrow k$ is defined as $\epsilon(T_\sigma )=1$ if $l(\sigma)$ is maximal (the maximal permutation has length $n(n-1)/2$ for the algebra on $n$ strands) and $\epsilon (T_\sigma)=0$ otherwise. With a slight abuse of notation, we  now refer to $\sigma$ or even $T_{\sigma}$ as maximal when $l(\sigma)$ is maximal.  The maximal $T_\sigma$ is called $max_n$.

\begin{prop} $\epsilon$ is a non-degenerate trace map and makes \NC{} and \Hn{} into Frobenius algebras.
\end{prop}
\begin{proof} $\epsilon$ non-degenerate means that for any $\sigma \in S_n$, there exist $\gamma, \gamma' \in S_n$ such that $T_\sigma T_\gamma =T_{\gamma'}T_\sigma=max_n$.  The properties of the length function on \Sn{} and the multiplication rules derived in  \ref{NC} and \ref{Hn} imply that this is indeed the case.

\end{proof}

The $0$-Hecke and Nilcoxeter algebras are not symmetric.  However, one can define an involution $f$ on $\A$ for both $\A=\NC$ and $\A=\Hn$: 
\begin{eqnarray*}
f(T_i)&=&T_{n-i}, 1\le i \le n-1\\
f(ab)&=&f(a)f(b).
\end{eqnarray*}

It is clear the $\forall \sigma$, $f(T_{\sigma})=T_{\gamma}$ for some $\gamma$, and that $l(\sigma)=l(\gamma)$.  It immediately follows that $\sigma$ is maximal if and only if $\gamma$ is maximal, so that $\epsilon(T_{\sigma})=0$ if and only if $\epsilon(T_{\gamma})=0$.  Hence, $\epsilon(a)=\epsilon(f(a))$ for all $a\in A$.
\\
\\

{\bf Acknowledgments:} I would like to thank my advisor Mikhail Khovanov for his guidance as well as the NSF for supporting me while writing this paper through the grant DMS-0706924.

\section{Topological Calculi of Diagrams}

In order to determine the dimension of the center of our Frobenius algebra $\A=\NC$ or $\A=\Hn$, we  make use of relation (1) of Claim \ref{TwistedDuality} and determine the dimension of \AT.  A basis element of \A{} can be represented as a monotonic immersion of $n$ unit intervals $[0,1]$ in $\real \times[0,1]$, subject to some relations.  We can see this element as $n$ strings each going from a  top position $1\le i\le n$ at $1$ to a bottom position $1\le j \le n$ at $0$.  A generator $T_i$ of \A{}  is represented by the crossing of the $i$ and $i+1$ strings.  The braid relations of \A{} allow us to change this diagram of strings using the third Reidemeister move ($T_iT_{i+1}T_i=T_{i+1}T_iT_{i+1}$) and to move far away crossings up and down with respect to each other ($T_iT_j=T_jT_i$ for $|i-j|>1$).  We read diagrams from top to bottom.

\begin{figure}[h]
\begin{center}
\includegraphics[width=5in]{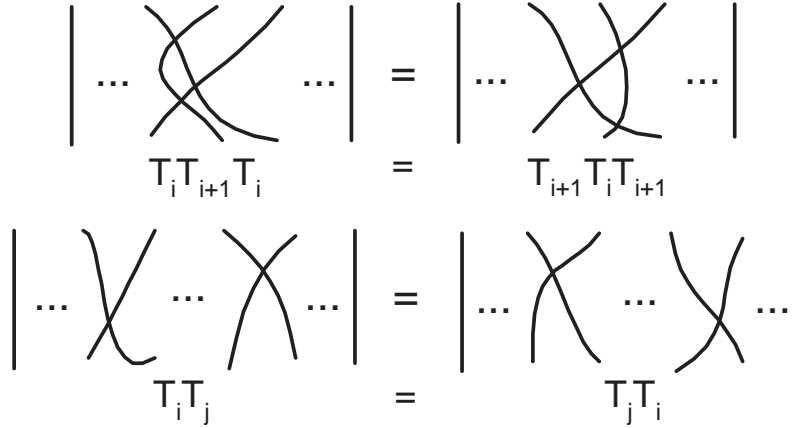}
\end{center}
\caption{These diagrams show the graphical relations commun to \C[\Sn], \NC{} and \Hn. The relation used to handle $T_i^2$ is what distinguishes these algebras.}
\label{Fig:EquivOnPlane}
\end{figure}

	\subsection{Diagrams on a Flat Band} 
A diagram like the one described above, a monotonic immersion of $n$ unit intervals $[0,1]$ in $\real \times[0,1]$, can be used to represent a basis element of the symmetric group ring \C[\Sn] if we impose the usual Reidemeister moves as relations.  By changing the relations appropriately, it can also represent a basis element of \NC{} or \Hn.  To $s_i\in \Sn$ or $T_i\in \A$ we assign a crossing of the $i$ and $i+1$ strands, to a product $s_is_j$ or $T_iT_j$ we assign a crossing of the  $i$ and $i+1$ strands followed by a crossing of the $j$ and $j+1$ strands and to a linear combination of basis elements we assign a formal sum of diagrams.  The relations imposed are:\\
for a basis element of \C[\Sn]:

\begin{enumerate}
\item
$s_is_{i+1}s_i=s_{i+1}s_is_{i+1}$
\item
$s_is_j=s_js_i$ for $|i-j|>1$
\item
$s_i^2=1$
\end{enumerate}

for a basis element of \NC:
\begin{enumerate}
\item
$T_iT_{i+1}T_i=T_{i+1}T_iT_{i+1}$
\item
$T_iT_j=T_jT_i$ for $|i-j|>1$
\item
$T_i^2=0$
\end{enumerate}

and for a basis element of \Hn:

\begin{enumerate}
\item
$T_iT_{i+1}T_i=T_{i+1}T_iT_{i+1}$
\item
$T_iT_j=T_jT_i$ for $|i-j|>1$
\item
$T_i^2=T_i.$
\end{enumerate}

\begin{figure}[h]
\begin{center}
\includegraphics[width=4in]{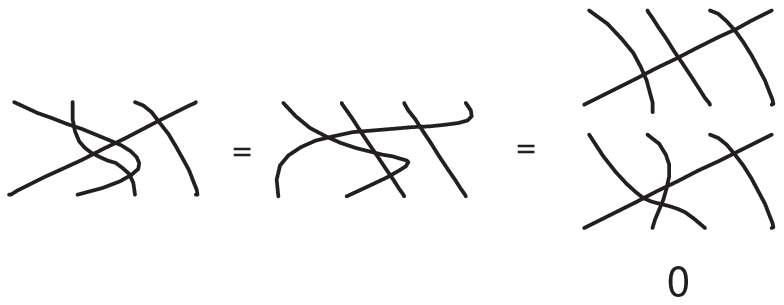}
\end{center}
\caption{This first diagram is equivalent to three distinct diagrams in \C[\Sn], \Hn{} and \NC{} respectively. }
\label{Fig:EquivComparison}
\end{figure}

\begin{figure}[h]
\begin{center}
\includegraphics[width=4in]{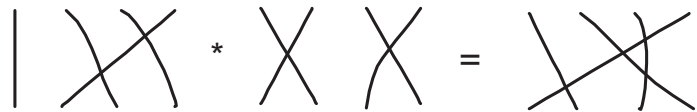}
\end{center}
\caption{Multiplication of basis elements corresponds to concatenation of diagrams.}
\label{Fig:Multiplication}
\end{figure}

Right multiplication of two basis elements $b_1,b_2$ then becomes the  juxtaposition of the diagram for $b_2$ under the diagram for $b_1,$
then rescaled to be in $\real \times[0,1]$ as shown in Figure \ref{Fig:Multiplication}.

To each word $w=s_{i_1}\dots s_{i_k}$ in \Sn{} there is assigned an element $T_w=T_{i_1}\dots T_{i_k}$ of \A{} and a diagram $D$.  In the Nilcoxeter algebra case, if $w$ is not reduced, then, $T_w=0$ and $D=0$.  In the $0$-Hecke algebra case, $T_w=T_{w'}$ for some reduced word $w'$ and likewise, $D=D'$ for some reduced diagram $D'$.  Note that in general $w$ and $w'$ are not equivalent words in \Sn.  Moreover, because a reduced word in \Sn{} can be related to any other equivalent reduced word through the braid relations alone, diagrams in \NC{} or \Hn{} corresponding to reduced words equivalent in \Sn{} are equivalent.  Conversely, if two reduced diagrams are equivalent, then the reduced words to which they correspond are also equivalent in \Sn.

	\subsection{Diagrams on a Cylinder: The Quotients of the Nilcoxeter and 0-Hecke Algebras by Their Commutator Subgroups}

\begin{figure}[h]
\begin{center}
\includegraphics[width=4in]{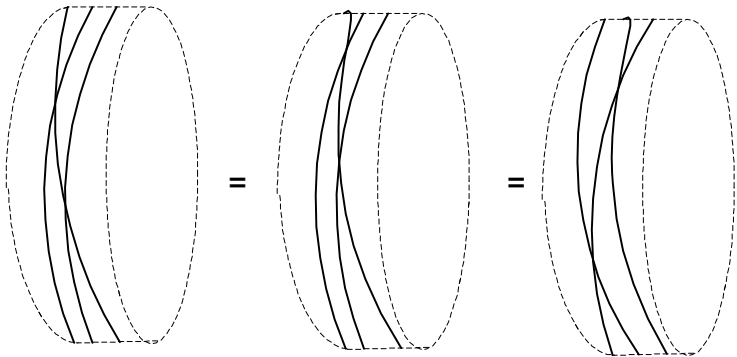}
\end{center}
\caption{Being on a cylinder imposes the relation $ab=ba$ since the top crossing can be carried around the cylinder to become the bottom crossing.}
\label{Fig:CylinderFullPic}
\end{figure}

Looking at diagrams on \real$\times S^1$ instead of \real $\times [0,1]$  imposes the new relation $ab=ba$ on \A{}  as shown in Figure \ref{Fig:CylinderFullPic}.  This is equivalent to working in \AC.

\begin{figure}[htb]
\begin{center}
\includegraphics[width=6in]{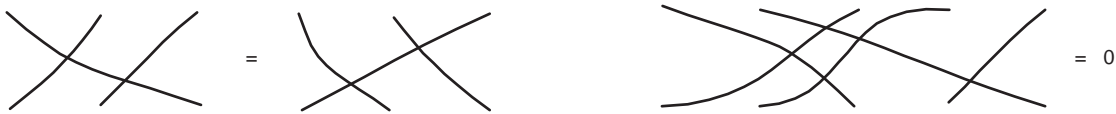}
\end{center}
\caption{These first two diagrams represent basis elements equivalent in \AC{} but not in \A{} for $\A=\NC$.  The third one represent the zero element in \AC{} for $\A=\NC$. This is because of the new relation $ab=ba$ which is imposed by working on a cylinder instead of a flat band.}
\label{Fig:OnCylinder}
\end{figure}

Figure \ref{Fig:OnCylinder} examplifies some of the consequences of this new relation for $\A=\NC$.  Recall that \AC{} is no longer an algebra but a vector space, both in the case $\A=\NC$ and $\A=\Hn$.

	\subsection{Diagrams on the \M: The Quotients of the Nilcoxeter and 0-Hecke Algebras by Their Twisted Commutator Subgroups}
		
\begin{figure}[h]
\begin{center}
\includegraphics[width=6in]{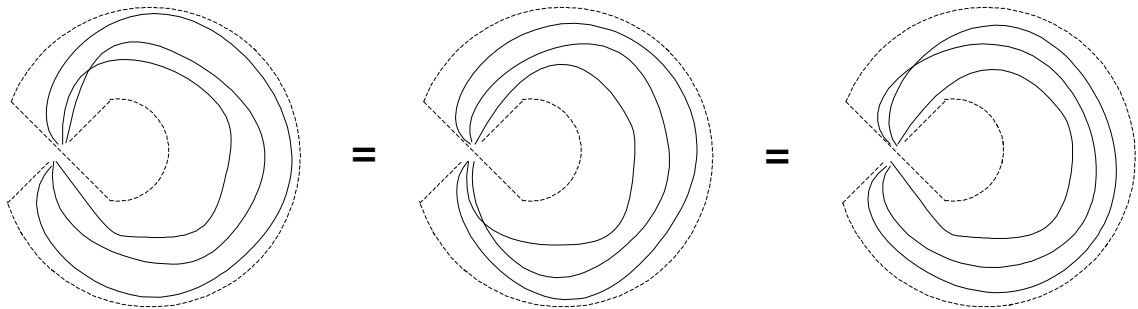}
\end{center}
\caption{Being on the \M $\,$imposes the new relation $ab=bf(a)$ with  $f(T_i)=T_{n-i}$ since a crossing $T_i$ can be dragged along the \M $\,$ to become $T_{n-i}. $}
\label{Fig:MoebiusFullPic}
\end{figure}

\begin{figure}[h]
\begin{center}
\includegraphics[width=4.5in]{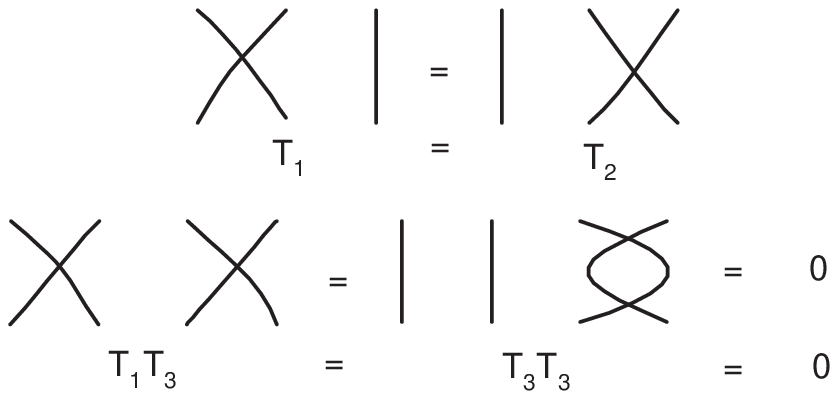}
\end{center}
\caption{The first two diagrams of this Figure are representations of the same element in \AT.  On the \M, if the left crossing is pushed once around, it becomes the right crossing.  This corresponds to the identification of diagrams in   \A{} under the  involution $f(T_i)=T_{n-i}$. The second diagram is 0 in \AT $\,$ for \A$=\NC$. }
\label{Fig:EquivOnM}
\end{figure}

 A basis element of \AT $\,$ can be represented as a diagram of $n$ strings living on a \M.  On the \M, a crossing $T_i$, if last in the diagram, can be pushed around the band to become a crossing $T_{n-i}$ first in the diagram. Figure \ref{Fig:MoebiusFullPic} shows how this can be done.  Thus, $TT_{i}=T_{n-i}T$ as desired in \AT. Indeed, on the \M, this new relation corresponds exactly to the involution used to define the relation on \AT. We  typically represent these diagrams on a flat band, but one should mentally identify the top and bottom of the band in opposite directions so as to form a \M.  Figure \ref{Fig:EquivOnM} gives an example of Nilcoxeter diagrams which are equivalent on the \M $\,$ but not on a flat band or on a cylinder. The elements they represent are therefore identified in \AT.  

	With the goal of  understanding the center of \A, we  classify these diagrams on the \M. The cases $\A=\NC$ and $\A=\Hn$ are very similar.  We  therefore carry out this classification in detail for $\A=\NC$ and later give a summary of the corresponding results for $\A=\Hn$.


\section{Nilcoxeter Diagrams on the \M} \label{NCOnM}

\begin{figure}[h!]
\begin{center}
\includegraphics[width=5in]{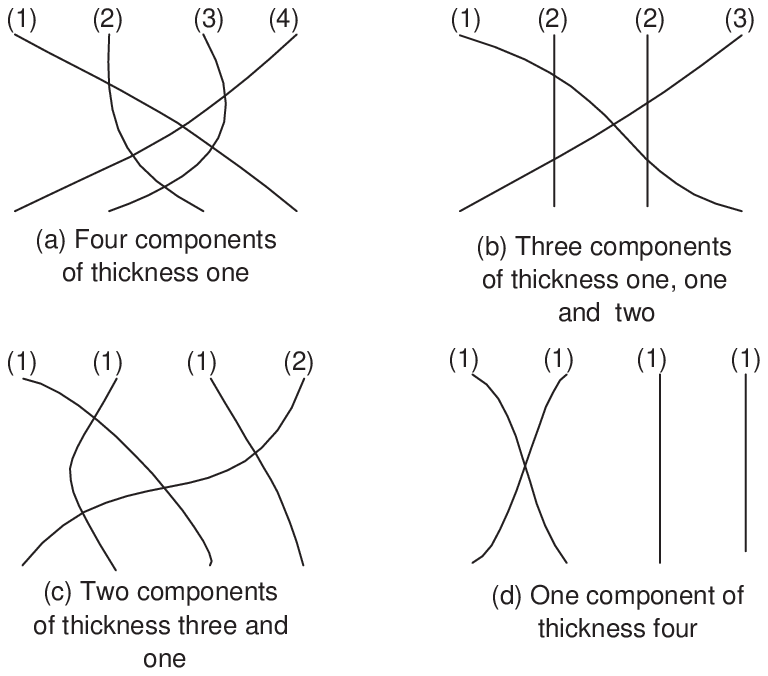}
\end{center}
\caption{These four non-zero diagrams are all of total thickness four but decompose into different prime components.}
\label{Fig:DiffNumberComps}
\end{figure}

To begin our study of diagrams on the \M, we  fix some notation. First note that on the \M,  as well as on the cylinder, these $n$ strings could actually be any number $c\le n$ of immersed circles.  Each circle is called a component of the diagram and the number of times a single immersed circle goes around the \M, i.e. the degree of that immersion,  is called the thickness of that component. The sum of the individual thicknesses of the components of a diagram is $n$, the thickness of the diagram. Figure \ref{Fig:DiffNumberComps} shows four diagrams of thickness four that decompose into components of different thicknesses.  In general, one can see that a diagram of thickness $2k$ with no crossing must consist of $k$ components each of thickness $2$; one of the components having both top and bottom positions 1 and $n$, another positions 2 and $n-1$, etc. The diagram corresponding to the image of $max_n$ in \AT, the diagram with the maximum number of crossings, must consist of $n$ components each of thickness $1$; one component having top position 1 and bottom position $n$, anothe top position 2 and bottom position $n-1$, etc. 

\begin{figure}[h!]
\begin{center}
\includegraphics[width=3in]{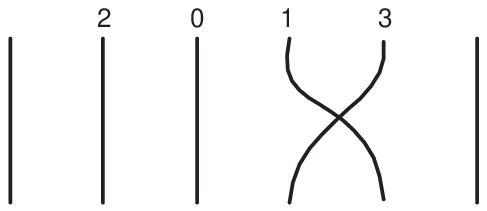}
\end{center}
\caption{The central component of thickness four of this diagram was here marked with the path $(p(0),p(1),p(2),p(3))=(3,4,2,5)$.  The other component has thickness two.}
\label{Fig:PathForComp}
\end{figure}

We now want to introduce the notion of a path $p$.  A path is defined in a specific component $C$ of a diagram $D$.  If $C$ has thickness $m$ and $D$ has thickness $n$, we define $p$ as a  sequence $(p(0), p(1),\dots ,p(m-1))$ where 
$1\le p(i)\le n,$ $p(i)\in\{\text{top position of a strand which belongs to }C\}$  is the top position of the strand of $C$ after going around the \M{} $i$ times starting at $p(0)$.  We  always consider this sequence $mod( n)$ and all statements made about it should be understood as such.  Figure \ref{Fig:PathForComp} gives an example of a path for the central component of thickness four of this composite diagram of total thickness six.  The starting position of this path is $p(0)=3.$  

For each component $C$, this sequence  is determined by $D$ once $p(0)$ is chosen.  Hence, two paths $p,p'$ in the same component can only differ by an integer $l\le m$, the lag of one path behind the other such that $p(i)=p'(i+l)$.

	\subsection{Prime Diagrams}

A prime diagram is a diagram consisting of a single component.    We will see that there is a unique prime diagram of thickness $n$ for each $n$.

		\subsubsection{The Minimum Number of Crossings of a Prime Diagram}

First, let us look at the minimal number of crossings such a diagram must have.  If, instead of being on the \M, we were on a simple cylinder $S^1\times[1,n]$, the minimal condition for the $n$ strands to be connected, to be part of the same component, would be clear: the diagram would have to correspond to an $n$-cycle in \Sn,  in particular, it would need to have at least the $n-1$ possible crossings $T_1$,\dots $T_{n-1}$.  In fact, it needs to have exactly each possible crossing once since on the cylinder, two strings  are part of the same component if and only if they cross exactly once.  Similarly, since the string with top position $i$ and the one with bottom position $n-i+1$ are already part of the same component in a diagram on the \M, one can quickly see that on the \M,  a diagram of thickness $n$ needs one of $T_i$ or $T_{n-i}$ for each $1\le i\le \lfloor (n-1)/2 \rfloor$ to be connected.  The minimum number of crossings required is therefore  $\lfloor (n-1)/2 \rfloor$.

		\subsubsection{The Exact Number of Crossings of a Prime Diagram}

We  need the following result to determine the maximum number of crossings that a non-zero prime diagram can have.

\begin{prop}{\em A diagram $D$ on the \M $\;$is zero if and only if there are more than one crossing between two fixed paths. }
\end{prop}
\begin{proof}
The first implication is clear: if there are two consecutive crossings $T_i$ for some $i$ in $D$, then these two crossings must occur between the same two paths.  
Conversely, given two paths $p$ and $p'$, one can see them as being both ``on top'' of the other strands of $D$ for each $i$.  One can then slide a crossing occuring between these two paths around the \M{}  to the next or previous $i$. The only obstacle to doing this is if there is already a crossing between these paths at this $i$, in which case we have two consecutive crossings and $D=0.$  If there are two crossings between these two paths occuring at different $i$, one can then move one of them around the \M{}    until they are consecutive and annihilate the diagram. Figure \ref{Fig:MoveCrossingOnPath} shows how this is done for a specific example. 

\begin{figure}[tbh]
\begin{center}
\includegraphics[width=5in]{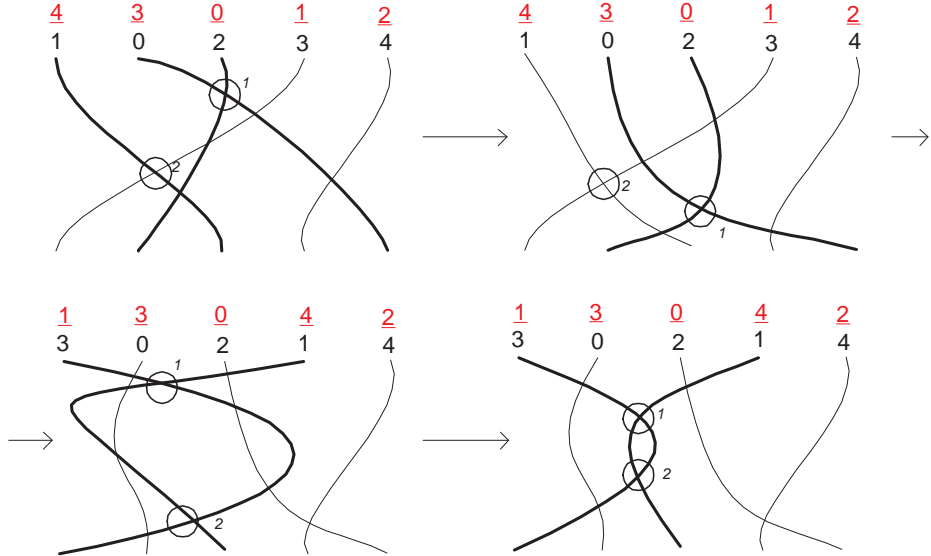}
\end{center}
\caption{Two paths and two crossings have been marked on this Figure.  Crossing 1 occurs at $i=0$ and crossing 2 at $i=1$. The bold strands are considered to be "on top" of the others, allowing us to see easily how to move crossings between them. One can first slide crossing 1 down the diagram. Then it can be pushed around the \M to $i=1$ so that the two crossings now both occur at $i=1$.  They become consecutive and $D$ is seen to be zero.  Note that these two paths actually cross a third time originally at $i=3$.}
\label{Fig:MoveCrossingOnPath}
\end{figure}
\end{proof}

We have already shown that in order for $D$ to be prime, it must have at least $\lfloor (n-1)/2 \rfloor$ crossings.  We now check that this is actually the maximum number of crossings a prime diagram of thickness $n$ can have if it is to be non-zero.

\begin{prop} \label{NumberCrossingsPrime}
A prime diagram of thickness $n$ has exactly  $\lfloor (n-1)/2 \rfloor$ crossings
\end{prop}
\begin{proof}
 Let us fix a component $C$.  We want to look at (unordered) pairs of paths in $C$, but some pairs are redundant.  For instance, in Figure \ref{Fig:PathForComp}, we would like the two pairs of paths $(p=(3,4,2,5),p'=(2,5,3,4))$ and $(q=(4,2,5,3),q'=(5,3,2,4))$ to be equivalent.  We  therefore consider pairs of paths modulo relative starting points, that is we impose the relation $$(p,p')\equiv(q,q')\Leftrightarrow  (p(i),p'(i))=(q(i+c),q'(i+c))$$ for some integer $c$ (remember that we are working $mod(m)$).  Note that if $p'$ lags on $p$ by $l$, then $p$ lags on $p'$ by $m-l$.  Since we are looking at unordered pairs and since any two pairs with the same lag are equivalent by the above relation, any pair of paths $(p,p')$ such that $p(i)=p'(i+l)$ is equivalent to any pair $(q,q')$ such that $q(i)=q'(i+n-l)$ since then $q'(i)=q(i+l)$.  The pairs of paths modulo relative starting points are therefore indexed by an integer $1\le l\le \lfloor m/2\rfloor$.  In our example of Figure \ref{Fig:PathForComp}, we thus have a total of $2$ classes of  pairs of paths, those which have a lag of $1$ or equivalently of $3$ such as $(p=(1,3,4,2),p'=(2,4,1,3))$, and those which have a lag of $2$ such as $(q=(1,3,4,2),q'=(3,1,2,4))$.\\
Remember that when the same pair of paths intersects twice, the diagram is zero.  Moreover, a crossing belongs to a pair of paths  if and only if it also belongs to all equivalent pairs. Thus, there can be only one crossing per equivalence class of pairs of paths 
for the component to be non-zero.  Since there are $\lfloor m/2 \rfloor$ equivalence classes of pairs of paths, there can be a maximum of $\lfloor m/2 \rfloor$ crossings within a component of thickness $m$ for the diagram $D$ to be non-zero. \newline
In this section, we are concerned with prime diagrams, so now let $D$=$C$.  For a prime diagram of odd thickness $n$, the maximum number of crossings we have just derived $\lfloor n/2 \rfloor=\lfloor (n-1)/2 \rfloor$, the minimum number of crossings previously derived, and the Proposition is proved.  For $D$ of even thickness $n$, we will show that $D$ cannot actually realize this maximum of $\lfloor n/2 \rfloor$ crossings; it has either no more than $n/2-1=\lfloor (n-1)/2 \rfloor$ crossings or is zero, which  completes the proof of Proposition \ref{NumberCrossingsPrime}. \newline
In fact, we show that a non-zero even prime diagram cannot have a middle crossing $T_{n/2}$.  First, note that the paths $p$ and $p'$, $p(0)<p'(0)$ never intersect iff 
$$p(2k)<p'(2k)$$ and 
$$p(2k+1)>p'(2k+1)$$
for all $k$ (again, we consider both 2k and 2k+1 as integers $mod(n)$). Now, assume that $D$ contains no middle crossing $T_{n/2}$ and consider the paths $(p,p'); p(0)=n/2, \; p'(0)=n/2+1$. Since no strand ever crosses the middle point of the diagram, we have 
$$p(2k)\le n/2, \; p'(2k)\ge n/2+1 $$
 and $$p(2k+1)\ge n/2+1, \; p'(2k+1)\le n/2$$
 such that $p(2k)<p'(2k)$ and $p(2k+1)>p'(2k+1)$ and $p$, $p'$ do not intersect. Having at least one class of pairs of paths which do not intersect, such a diagram $D$ has a maximum of $n/2-1=\lfloor (n-1)/2 \rfloor$ crossings.  Assume then that $D$ does contain a middle crossing $T_{n/2}$.  Without loss of generality, assume that $T_{n/2}$ is the top crossing of the diagram and consider the same paths $(p,p'); p(0)=n/2, p'(0)=n/2+1$. Now, since these paths cross immediately, $p(1)\le n/2<n/2+1\le p'(1)$.  Suppose that these paths do not intersect again, then 
$$p(2k)>p'(2k)$$ 
and $$p(2k-1)<p'(2k-1)$$
 for $ k>0$.  But we know that 
$$n/2=p(0)=p(n)=p(2(n/2))<p'(2(n/2))=p'(n)=p'(0)=n/2+1,$$ 
which is a contradiction.  Therefore, (p,p') must intersect more than once and $D$ is zero.
\end{proof}

		\subsubsection{The Crossings of a Prime Diagram}

We have now shown that a prime diagram of thickness $n$ must have exactly $\lfloor (n-1)/2 \rfloor $ crossings.  Moreover, we know that these crossings are exactly  $T_i$ or $T_{n-i}$ for each $1\le i\le  \lfloor  (n-1)/2 \rfloor$.  In any case, each $T_i$ appears at most once, which  allows us to push all of the crossings to the left side of the diagram and thus assume that the $\lfloor n/2 \rfloor$ crossings are exactly $T_i$ for $1\le i\le \lfloor (n-1)/2 \rfloor$.  Indeed, to see this first assume inductively that all $T_j$, for $ j>n-i$, have already been pushed to the left side. We can do this since there is at most one $T_{n-2}$ crossing that could prevent $T_{n-1}$ from  being moved. If there is a $T_{n-2}$ crossing above or below $T_{n-1}$, the latter is free to be pushed in the other direction and if there is no $T_{n-2}$ then the $T_{n-1}$ crossing can be pushed either up or down the \M{} to become $T_1$. Consider now what could block $T_{n-i}$.  Everything to its right has already been pushed, so only  $T_{n-i-1}$ above or below it can block it. It is therefore free to be pushed in the other direction. Therefore, from now on when we think of a prime diagram, we  always represent it with its crossings being $T_1, T_2,\dots ,T_{\lfloor (n-1)/2 \rfloor}$ in some order.  It remains to show that these $\lfloor (n-1)/2\rfloor!$ diagrams are in fact all equal in \AT.
			
		\subsubsection{Uniqueness of Prime Diagrams}

\begin{lem} {\em $T_1T_2\dots T_{\lfloor (n-1)/2 \rfloor}=T_{\sigma(1)}T_{\sigma(2)}\dots T_{\sigma(\lfloor (n-1)/2 \rfloor)}\in$ \AT{}  for all $ \sigma \in S_{\lfloor (n-1)/2 \rfloor}$}
\end{lem}
\begin{proof}
First, let us look at the case when $\sigma$ is a cyclic permutation.  Fix $m={\lfloor (n-1)/2 \rfloor}$. Now,
$${\bf T_1}(T_2\dots T_{m})=(T_2\dots T_m){\bf T_{n-1}}={\bf T_{n-1}}(T_2\dots T_m)=(T_2\dots T_m){\bf T_{n-(n-1)=1}}$$
since $|(n-1)-i|>1$ for $2\le i\le m$.
More generally, if we have any diagram $D=T_{\sigma'(1)}T_{\sigma'(2)}\dots T_{\sigma'(m)}$ for $\sigma'\in S_m$, then
\begin{align*}
{\bf T_{\sigma'(1)}}(T_{\sigma'(2)}\dots T_{\sigma'(m)})=&(T_{\sigma'(2)}\dots T_{\sigma'(m)}){\bf T_{n-\sigma'(1)}}=\\
{\bf T_{n-\sigma'(1)}}(T_{\sigma'(2)}\dots T_{\sigma'(m)})&=(T_{\sigma'(2)}\dots T_{\sigma'(m)}){\bf T_{\sigma'(1)}}\;\;\;\; (*) 
\end{align*}
because in fact $n-\sigma'(1)>m+1$ so that $|(n-\sigma'(1))-i|>1$ for $1\le i\le m$ unless $n$ is odd and $\sigma'(1)=m$.  Then $n-\sigma'(1)=m+1$, but we still have
that $|(m+1)-\sigma'(j)|>1$ for $2\le j\le m$ since $\sigma'(j)\ne m$ for $j\ne 1$ so that (*) still holds.  Diagramatically, we have taken the top crossing and pushed it onto the right side of the diagram, where it is free to be pushed to the top again.  We can then push it back to the bottom left side of the diagram.  Hence the Lemma holds for $\sigma$ a cyclic permutation.  Figure \ref{Fig:CyclicPerm} shows how this is done for a simple example.\newline

\begin{figure}[h]
\begin{center}
\includegraphics[width=4.5in]{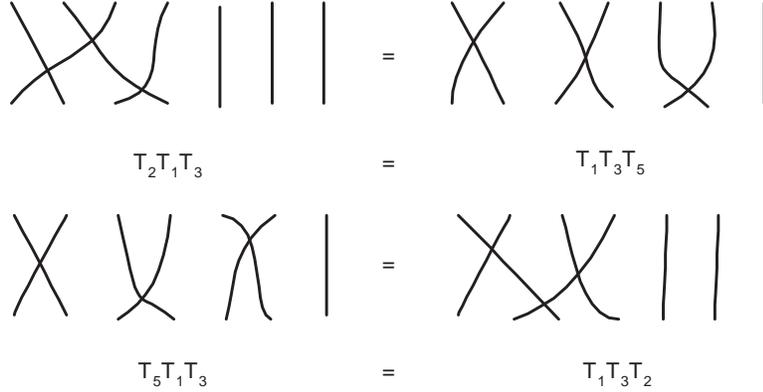}
\end{center}
\caption{This Figure shows how one can cyclically permute the crossings of a non-zero prime diagram. }
\label{Fig:CyclicPerm}
\end{figure}

Now let us compose a  permutation $\sigma'$ with a transposition $(i\;\;i+1)$. Since we already know that this diagram is invariant under cyclic permutations of its crossings, we can actually assume that $i=1$.  If $\sigma'(2)\ne\sigma'(1) \pm 1$, $\sigma'(1) $ and $\sigma'(2)$ commute and we are done.  Without loss of generality, assume therefore that $\sigma'(2)=\sigma'(1)+1$.  Let $k$ be such that $\sigma'(k)= \sigma'(2)+1$. So 
\begin{align*}
{\bf T_{\sigma'(2)}T_{\sigma'(1)}}(\dots T_{\sigma'(k)}\dots T_{\sigma'(m)})=\\(T_{n-\sigma'(k)}\dots T_{n-\sigma'(m)})({\bf T_{\sigma'(2)}T_{\sigma'(1)}})(\dots T_{\sigma'(k-1)})=\\
{\bf T_{\sigma'(2)}}(T_{n-\sigma'(k)}\dots T_{n-\sigma'(m)}{\bf T_{\sigma'(1)}}\dots T_{\sigma'(k-1)})=\\(T_{n-\sigma'(k)}\dots T_{n-\sigma'(m)}{\bf T_{\sigma'(1)}}\dots T_{\sigma'(k-1)}){\bf T_{\sigma'(2)}}=\\
(T_{n-\sigma'(k)}\dots T_{n-\sigma'(m)}){\bf T_{\sigma'(1)}T_{\sigma'(2)}}(\dots T_{\sigma'(k-1)})=\\{\bf T_{\sigma'(1)}T_{\sigma'(2)}}(\dots T_{\sigma'(k)}\dots T_{\sigma'(m)})
\end{align*}
using the fact that it is invariant under cyclic permutations.  Thus this diagram is in fact invariant under all permutations $\sigma \in S_{\lfloor (n-1)/2 \rfloor}.$  
\end{proof}


This Lemma now allows us to conclude that there is exactly one prime diagram of thickness $n$ for each $n$.  This diagram can be represented as $D=T_1T_2\dots T_{\lfloor (n-1)/2 \rfloor)}$.
			
	\subsection{Composite Diagrams}

We now need to see how prime diagrams can be combined into non-zero composite diagrams.  At first glance, there seems to be a number of ways to combine diagrams.  Observe that, given two diagram $D$ and $D'$ with the same prime components $C_1,\dots ,C_k$, it is always possible to go from $D$ to $D'$ via an isotopy. Indeed, given that there is no under and overcrossings, one can first consider that $C_1$ is on top and place it where one wants, then do the same with $C_2$ and thus simply rearrange the components of $D$ as they are in $D'$. This isotopy, however, might very well go through zero in  \AT{}  if $\A=\NC$.  In fact, we know that if the second Reidemeister move is used, then it does go through zero, and we also know that this move is the only one which changes the number of crossings.    Therefore, given isotopic diagrams $D$ and $D'$,  if $D$ has more crossings than $D'$, then $D=0$.  Thus all non-zero diagrams in \AT{} with the same prime components $C_1,\dots ,C_k$ have the same number of crossings and that number is the minimum of the isotopy class.  Note that such a minimum must exist and that the corresponding diagram, having no two paths crossing more than once since it is minimal, is non-zero.

Note also that the relations on these diagrams being local, any relation on a prime component continue to hold when the component is combined in a composite diagram.  
		
		\subsubsection{Combining Even Prime Diagrams}
	
\begin{figure}[h]
\begin{center}
\includegraphics[width=5in]{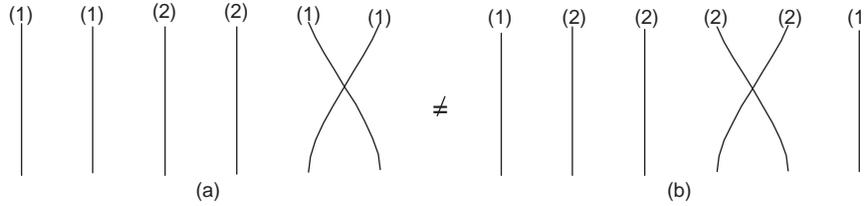}
\end{center}
\caption{These two diagrams are both non-zero since they are minimal and their decomposition into prime components is the same.  However, they do not represent the same element of \AT.  This shows that prime component decompositions does not quite determine unique elements in \AT. }
\label{Fig:EvenMincomp}
\end{figure}

Consider a diagram $D$ with even prime components $C_1, C_2$ of respective thicknesses $2k_1, 2k_2$.  Since even prime diagrams have no middle crossings, we can combine $C_1$ and $C_2$ by either nesting $C_2$ inside $C_1$, which means having $C_1$ at positions $1,2\dots ,k_1,k_1+2k_2+1,k_1+2k_2+2\dots ,2k_1+2k_2$  both on top and bottom and $C_2$ at positions $k_1+1,k_1+2,\dots k_1+2k_2$ or vice versa.  Neither of these two arrangements will add crossings between $C_1$ and $C_2$: if $C_1$ has $k_1-1$ crossings and $C_2$ has $k_2-1$ crossings and $D$ has $k_1+k_2-2$ crossings in either arrangement, which are therefore both minimal.  Let $D$ be the diagram with $C_1$ in the center and $D'$ the one with $C_2$ in the center.  Then an isotopy from $D$ to $D'$ requires the second Reidemeister moves which creates double crossings and the isotopy therefore goes through zero.  Hence, $D$ and $D'$ are distinct non-zero diagrams in \AT.	
	
\begin{figure}[h]
\begin{center}
\includegraphics[width=3in]{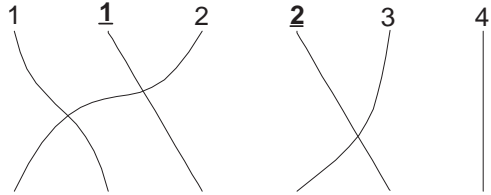}
\end{center}
\caption{This is a non minimal arrangement of two even prime components.  The strands are numbered with their positions in their respective prime components to show their assignment functions $f_{C_i}$ more clearly.  The inter-component crossings can be avoided by splitting one of the components in two and placing the other in the center.  This diagram, not being minimal, is zero.}
\label{Fig:NonMinEvencomp}
\end{figure}

		\subsubsection{Combining Odd and Even Components}

\begin{figure}[h]
\begin{center}
\includegraphics[width=3in]{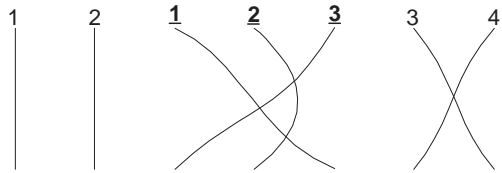}
\end{center}
\caption{This is a minimal diagram with one odd (composite) and one even prime component.  The strands are numbered with their positions in their respective prime components to show their assignment functions $f_{C_i}$ more clearly. Minimal arrangements of this kind always involve placing the odd component in the center so as to not create inter-component crossings.}
\label{Fig:MinEvenOdd}
\end{figure}

Next, given a component $K$ (prime or composite) of total odd thickness with $k$ crossings, and  even prime components $C_1,\dots ,C_m$, with $n_1,\dots ,n_m$ crossings respectively, the minimal arragements have no intercomponent crossings: they have $k+n_1+\dots +n_m$ crossings.  These arrangements all involve $K$ in the center, disjoint from the even components, which are split each in the middle around $K$.  If $C_1,\dots ,C_m$ are all distinct, then there are $m!$ arrangements of this type.    By the same argument as in the case of two even diagrams, there is no isotopy between these arrangements that does not require the second Reidemeister move, so these $m!$ diagrams are in fact distinct. Figure \ref{Fig:MinEvenOdd} is an example of such a minimal arrangement.
Since  $K$ has odd thickness and thus cannot be split in its middle, if it is not in the center, then it is also not disjoint from the even components. The corresponding diagram is then not minimal and is therefore zero.  
			
		\subsubsection{Combining Odd Prime Diagrams}
			
Let us now see how to combine two odd diagrams $C_1,C_2$ of thicknesses $2k_1+1,2k_2+1$ respectively into a diagram $D$.  To try to find a combination, one can start by numbering $2k_1+1$ of the $2k_1+2k_2+2:=n$ top positions of $D$ with numbers $1,..,2k_1+1$ in a certain color.  These would be the top position for the strands of $C_1$.  First note that once these top positions for $C_1$ are chosen, the bottom ones are determined.  Indeed, if the number $i$ was assigned to the top position $j$ of $D$, then the number $n-i$ is assigned to the bottom position $n-j+1$.  We call this assignment from $\{1,2,\dots ,2k_i+1\} \hookrightarrow  \{1,2,\dots ,n\}$, $f_{C_i}.$  There seems to be many ways one can choose the positions of $C_1$ and $C_2$ in $D$, but we will see that there is in fact a unique choice that makes $D$ non-zero. 

\begin{lem} {\em There is a unique non-zero diagram $D$ with prime components $C_1,\dots C_m$ of respective thicknesses $2k_1+1,\dots 2k_m+1$.}
\end{lem}
\begin{proof}  We  describe an arrangement and then argue that it is the unique one with the minimum number of crossings.  Without loss of generality, assume that $$2k_1+1\le 2k_2+1\le \dots \le2k_m+1.$$  Let $n:=2k_1+2k_2+\dots +2k_m+m$ be the total thickness of $D$.  We  build the diagram component by component.  Let $D_i$ be the diagram of thickness $n_i :=2k_1+2k_2+\dots +2k_i+i$ with prime components $C_1,\dots ,C_i$ obtained after the $i^{th}$ step of this procedure.  The following assignment functions are the ones yielding the unique minimal diagram $D_{i+1}$:
$$f_{C_{i+1}}(j)=\left\{\begin{array}{ll}j&j\le  k_{i+1} \\j+n_i+k_{i+1}&j>k_{i+1}\end{array}\right.$$ and
$$f_{D_i}(l)=l+k_{i+1},$$
which means that we keep the previous composite $D_i$ in the center and give $C_{i+1}$ the edge positions.

\begin{figure}[htb]
\begin{center}
\includegraphics[width=5in]{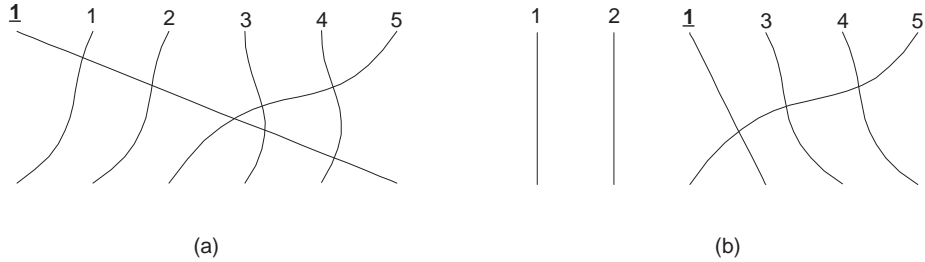}
\end{center}
\caption{The strands of these two diagrams are numbered with their positions in their respective prime components to show their assignment functions $f_{C_i}$ more clearly. Diagram (a) is non-minimal because the thickest odd component was placed in the center of the diagram, creating $1\times5$ inter-component crossings, while diagram (b), by placing the thickest component on the outside, minimized inter-component crossings to $1\times 1$. }
\label{Fig:MinAndNonOdd}
\end{figure}

First, we want to minimize the number of extra crossings of a component with itself.  If possible, it is clear that a minimal diagram is one where the placement of each component does not add crossings between strands of the same component.  If $f$ is not monotonic, extra crossings might  be created, while none 
are added if $f$ is monotonic. Therefore, monotonic assignment functions create minimal diagrams.  Note that for each component of thickness $m$, if $f'=m+1-f$ then $f$ and $f'$ are equivalent on the \M.  Indeed,
 diagrams which are reflections of each other are seen to be the same by rotating the \M. We  therefore only consider strictly increasing assignment 
functions $f$.  Moreover, all assignment functions $f_{C_i}$ for $C_i$ a prime component which do not increase the number of crossings within $C_i$ itself are equivalent since all diagrams of $C_i$ with the minimal number of crossings are equivalent.

\begin{figure}[h!]
\begin{center}
\includegraphics[width=2in]{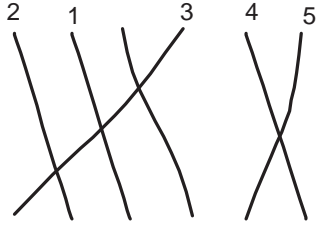}
\end{center}
\caption{The assignment function $f_{C_i}$ shown by the numbered strands is not monotonic and does not preserve the number of crossings of that prime component of thickness five. This diagram is therefore zero. }
\label{Fig:NonMonotonic}
\end{figure}

We also want to minimize crossings between $D_i$ and $C_{i+1}$, which means that we want to avoid alternating positions,ie. situations 
like $$f_{D_i}(j)<f_{C_{i+1}}(k)<f_{D_i}(l)$$ which  add intercomponent crossings.  
For a component $C$ of thickness $m$ in a diagram of total thickness 
$n$, the botttom $j$ position is the same as the top $j$ position whenever $$f(j)= n+1-f(m+1-j).$$  When this is not true, we have additional 
intercomponent crossings.  Therefore, minimizing the crossings means that we want one component in the center and one on the outside of the diagram, 
each symmetric in their assignment. 

\begin{figure}[h!]
\begin{center}
\includegraphics[width=2in]{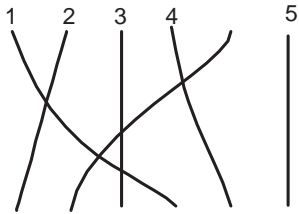}
\end{center}
\caption{The assignment functions $f_{C_i}$ shown by the numbered strands is not symmetric and does not minimize the number of intercomponent crossings. This diagram is therefore zero. }
\label{Fig:NonSymmPos}
\end{figure}
 
Inductively, we see that in any diagram, there is one strand per odd component crossing the center line of the 
diagram.  Therefore, if  $D_i$ is placed as the outside component in $D_{i+1}$, there are $i\times (2k_{i+1}+1)$ crossings between 
$D_i$ and $C_{i+1}$, since $D_i$ is composed of $i$ odd prime components and one strand for each of these prime components will cross the $2k_{i+1}+1$ strands of $C_i$.  Now, an odd component cannot really be split in two and placed on the outside because of the symmetry requirement.  Placing it on 
the outside therefore implies taking the first and last $k$ position, both top and bottom to minimize crossings.  The remaining position  has bottom 
position assignment $n_{i+1}+1-j$ if it has top assignment $j$.   $C_{i+1}$ being an odd prime component, we know that, in the standard 
representation of this diagram,  the strand having top position $k_{i+1}+1$ has bottom position  $1$ in $C_{i+1}$ and the one having bottom position 
$k_{i+1}+1$ has top position $k_{i+1}$.  This means that in $D_{i+1}$, if $f_{C{i+1}}(k_{i+1})=j$, the strand having top position $k_{i+1}+1$ in $C_{i+1} $ will cross $D_i$ $j-k_{i+1}$ times and the one with bottom position $k_{i+1}+1$ in $C_{i+1}$ will cross $D_i$ $n_{i+1}+1-j-(k_{i+1}+1)$ times for a total  of $$j-k_{i+1}+n_{i+1}+1-j-(k_{i+1}+1)=n_{i+1}-(2k_{i+1}+1)=n_i$$ times.  Now, $i\times (2k_{i+1}+1)>n_i$ since we have ordered the components such that 
$$k_1\le k_2\le\dots .$$  Therefore, placing the (thickest) prime odd component $C_{i+1}$ as the outside component always yields a minimal diagram.  It remains to 
see whether the diagram is invariant under  the choice of $f_{C{i+1}}(k_{i+1})=j$, the position in $D_{i+1}$ of the middle position of $C_{i+1}$.  This can be seen by pushing the crossing $T_{j-1}$ around the \M{} so 
that it becomes $T_{n_i+1}+1-{j-1}$.  The diagram thus becomes exactly the diagram with $f_{C{i+1}}(k_{i+1})=j-1$.  Doing the opposite, we get the 
diagram with $f_{C{i+1}}(k_{i+1})=j+1$, and so we see that $D_{i+1} $ is invariant under the choice of $j$.  There is therefore a unique diagram $D$ on the \M{} with prime components $C_1,\dots ,C_m$.  
\end{proof}

This completes the classification of Nilcoxeter diagrams on the \M.  We have learned that there is a unique prime diagram of thickness $n$ for each $n$, that there is a unique diagram with only odd prime components $C_1,\dots,C_m$ and that in general, the number of arrangements depends on the number of  even prime components of each thickness. 

\section{The Center of the Nilcoxeter Algebra}
	\subsection{The Dimension of the Center}
		We have now enumerated all the possible arrangements of prime components into composite diagrams.  Given even prime components $C_1,\dots, C_m$, $i_j$ of them with thickness $2j$  and odd prime components $C'_1,\dots ,C'_m$, one can make exactly  $\frac{m!}{\Pi_j (i_j!)}$ distinct non-zero diagrams $D$ with this prime decomposition.  We can now count the number of non-zero diagrams of thickness $n$ in \AT, each of which is a basis element of \AT.  This number is also the dimension of the center $Z$ of $\A$:
$$dim(Z)=\sum_{\lambda \vdash n}\frac{ n_{\lambda}!}{m_{\lambda}} ,$$
where $n_{\lambda}$ is the number of even parts in $\lambda$ and $m_{\lambda}=\Pi_j( i_j!).$

For instance, for $n=3$ we have partitions $$(1,1,1), (1,2), (3)$$ and the formula gives us
$$dim(Z)=0!/0!+1!/1!+0!/0!=3.$$
A more interesting case is $n=6$ where we have partitions 
\begin{eqnarray*}
(1,1,1,1,1,1), (1,1,1,1,2), (1,1,1,3), (1,1,2,2), (1,1,4), (1,5), (1,2,3),(2,2,2), (2,4), (3,3),(6).
\end{eqnarray*}
  The formula then gives us
$$dim(Z)=$$ 
$$0!/0!+1!/1!+0!/0!+2!/2!=1!/1!+0!/0!+1!/1!+3!/3!+2!/1!1!+0!/0!+1!/1!=12.$$
Here the partition $(2,4)$ gives us two distinct diagrams.

	\subsection{A Basis for the Center}
		We now know the dimension of the center of \NC, but we would also like to have a basis for this center.  Recall that we used the fact that 
$$ Z(\A)\backsimeq  \AT^*$$
and that this isomorphism is given by the trace map $\epsilon$.  Therefore, we just need to find the dual of the diagrams we counted in the previous section.
First, given an element $T_{\alpha}\in \A$, we call $T_{\beta}$ complementary to $T_{\alpha}$ if $\epsilon(T_{\alpha}T_{\beta})=1$, $\alpha, \beta \in \Sn$.  Note that for $\A=\NC,$ this is equivalent to $\alpha\beta=max_n,$ which also means that there is a unique complementary $T_{\beta}$ for each $T_{\alpha}.$

\begin{figure}[h!]
\begin{center}
\includegraphics[width=5in]{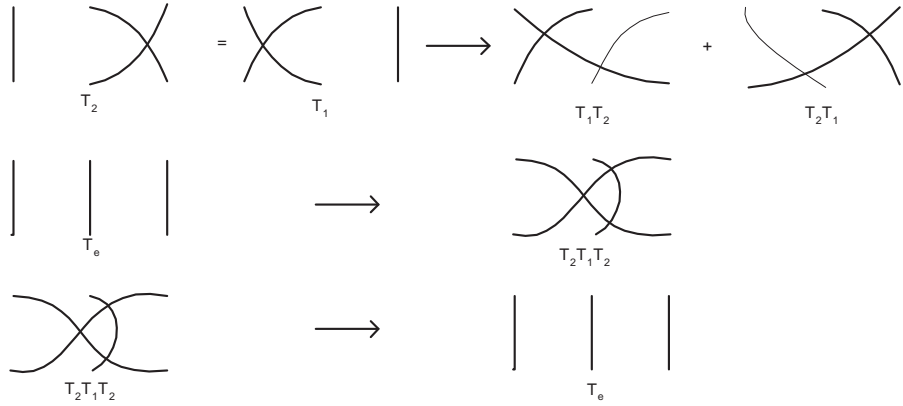}
\end{center}
\caption{ This figures shows the cosets corresponding to diagrams on the \M{} of total thickness three and the basis elements of $Z(\NC)$ formed by taking the sum of their complements. }
\label{Fig:DualNilc3}
\end{figure}

Now, we need to find the element of $Z(\A)$ corresponding to a basis element of \AT. Since the distinct diagrams of thickness $n$ formed a basis for the vector space \AT, their corresponding elements in $Z(\A)$ are also linearly independent and therefore form a linear basis for $Z(\A)$. The element $z$ of $Z(\A)$ corresponding to a given basis element $\overline{b}$ is such that $\epsilon(bz)=1$ for any $b\in \A$ with $\pi(b)=\overline{b}$, $\pi:\A\to\AT$ being the projection map.  To find these elements of $Z(\A)$, we first need to find the elements of the coset of \A{} corresponding to a given  element of \AT.  The element of $Z(\A)$ we want is then the sum of the complementary elements of the elements of that coset.  Here is an example for $n=3$.   Take the identity element in \AT{} represented by three non-intersecting strands on the \M.  Its coset in \A{} is just itself: $T_e$.  Its complementary element is $T_{max_n}$, our first basis element for $Z(\A)$.   Conversely, the coset of $T_{max_n}\in \AT$ is again just $T_{max_n}\in \A$ and its complementary element is the identity $T_e\in Z(\A)$, our second basis element.  Now let us consider $T_1\in \AT$.  Its coset is $\{T_1,T_2\}\subset\A$ since $T_1$ is equivalent to $T_2$ in \AT.  Their complementary elements are $T_{s_2s_1}=T_2T_1$ and $T_{s_1s_2}=T_1T_2$ respectively.  Our third basis element is therefore $T_2T_1+T_1T_2$.  This ends our computation of the basis of $Z(\A)$ for $n=3$.  See Figure \ref{Fig:DualNilc3} for the diagrams corresponding to this computation.

The Nilcoxeter algebra is a graded algebra with the grading given by the length function on \Sn.  Equivalently, the grading of a basis element is given by the number of crosssings in its corresponding diagram.  Note that the process described above to find a basis for $Z(\A)$  yields a homogeneous basis since every element of the coset of an element of \AT{} has the same length, or the same number of crossings.

	\subsection{Multiplication in the Center of \NC }
	
\begin{figure}[h]
\begin{center}
\includegraphics[width=3in]{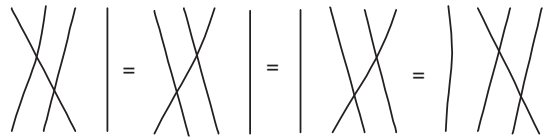}
\end{center}
\caption{ These are the diagrams in the coset corresponding to the partition $4=3+1$.  Note that the prime component of thickness three is on the outside of the odd component of thickness one. Since the prime diagram of thickness three always has a vertical strand either as its first or last strand, so do all these composite diagrams.}
\label{Fig:31Coset}
\end{figure}

We would now like to know what the multiplication table looks like in the basis we have described in the previous section.\\   
\begin{lem}{\em Every diagram corresponding to a basis element of \AT{} other than $T_{max_n}$ has a vertical line as its first or last strand.  Equivalently, a basis element of \AT{} other than $T_{max_n}$ cannot have both $T_1$ and $T_{n-1}$ in its expression.}
\end{lem}
\begin{proof}
First, let us look at a non-zero diagram $D$ on the \M{} corresponding to a partition of $n$ which has at least one even part.  Then, we know that one of these even prime diagrams is on the outside of $D$ and disjoint from its other components.  We also know that all the elements in the coset of this prime even diagram have at least one such vertical strand as either their first or last strand.  Since this even component is disjoint from all the other components of the diagram, this vertical strand clearly remains in all elements of the coset of $D$.  Figure \ref{Fig:411Coset} shows this for all the diagrams in the coset corresponding to the partition $6=4+1+1.$

\begin{figure}[h]
\begin{center}
\includegraphics[width=3in]{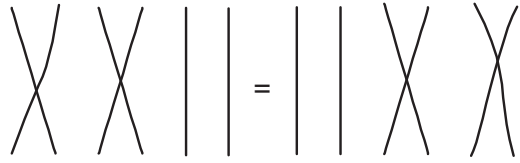}
\end{center}
\caption{These are the diagrams in the coset corresponding to the partition $6=4+1+1$.  Note that the even component of thickness four is always disjoint from the odd components of thickness one. Since the prime diagram of thickness four always has a vertical strand either as its first or last strand, so do all these composite diagrams.}
\label{Fig:411Coset}
\end{figure}

We now look at a non-zero diagram $D$ on the \M{} which is composed solely of odd prime components.  Since $D$ is not the diagram for $T_{max_n},$ one of its prime components has thickness greater than one. Let the thickest prime component $C$ of $D$ have thickness $2k+1.$ Remember from the classification of these diagram that this means that this component $C$ must be the component occupying positions $1$ and $n$ both at the top and at the bottom of the diagram. This is true for any diagram equivalent to $D$.   Once again, every diagram of an odd prime component must have a vertical strand at position 1 or $2k+1$. This is true also of $C$ in $D$ since $C$ occupies positions 1 and $n$ both at the top and bottom of the diagram.  Figure \ref{Fig:31Coset} shows that this is true in the case of the diagrams corresponding to the partition $4=3+1$.
\end{proof}

\begin{prop}{\em The multiplication table for the center of \NC{} in the basis described above is trivial.}
\end{prop}

\begin{proof}
 This vertical strand in every element of the cosets of our basis elements means that their complement must have either a strand going from top position $1$ to bottom position $n$ or from top position $n$ to bottom position $1$, depending on whether the vertical strand is at position $1$ or $n$, since $T_{max_n}$ has both of these diagonals. Since every basis element of $Z(\A)$ is a sum of such complements, every diagram appearing in any basis element of $Z(\A)$ has such a diagonal.

\begin{figure}[h]
\begin{center}
\includegraphics[width=3in]{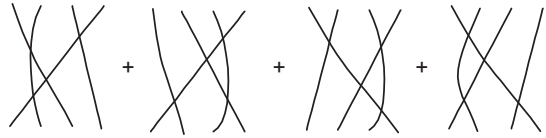}
\end{center}
\caption{ This is the basis element of $Z(\NC)$ for $n=4$ corresponding to the partition $4=3+1$.  It is the sum of the complements of the elements in the corresponding coset in \NC. }
\label{Fig:31Complement}
\end{figure}

Multiplying two elements with the opposite diagonals obviously yields a non-minimal, and therefore zero, diagram since the strand in the product corresponding to these diagonals  crosses every other strand twice.  

\begin{figure}[h]
\begin{center}
\includegraphics[width=2.5in]{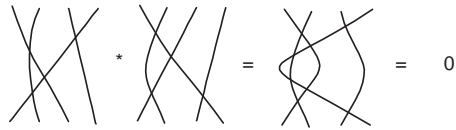}
\end{center}
\caption{ This is an example of the multiplication of two elements of \NC{} with opposite full diagonals.  Such a product is always zero. }
\label{Fig:MultOppDiag}
\end{figure}

If we multiply two elements with the same diagonal, we also have a trivial diagram.  Indeed, remember that if two strands cross more than once, the diagram is zero.  Since this diagonal strand crosses all other strands in the first element and does not remain vertical in the second because of the full diagonal of this second diagram, it must cross at least one strand twice.

\begin{figure}[h!]
\begin{center}
\includegraphics[width=2.5in]{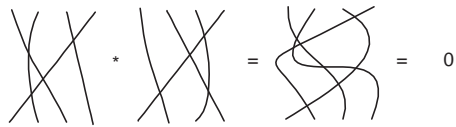}
\end{center}

\caption{ This is an example of the multiplication of two elements of \NC{} with the same full diagonal.  Such a product is always zero since the strand corresponding to the diagonal in the first element crosses at least one strand twice. }

\label{Fig:MultSameDiag}
\end{figure}

Hence, the product of any basis element by any element other than the identity is zero and the multiplication table of $Z(\A)$ is trivial as claimed.

\end{proof}

\section{The Center of the 0-Hecke Algebra}
	\subsection{0-Hecke Diagrams on the \M }
		
The classification of \Hn{} diagrams on the \M{} band is identical to that of \NC{} diagrams.  All of the arguments of Section \ref{NCOnM} go through for \Hn.  Everytime a diagram was eliminated for being zero in \NC, it is eliminated for not being minimal in \Hn.  Zero diagrams in \NC{} correspond exactly to non-minimal diagrams in \Hn.  Non-zero diagrams which are equivalent in \NC{} are also equivalent in \Hn{} and non-minimal diagrams which are equivalent in \Hn{} are equivalent in \NC.  Hence, the vector spaces \AT{} for $\A=\NC$ and $\A=\Hn$ are isomorphic and their bases are represented by the same diagrams on the \M.

Therefore, the dimensions of the centers of \NC{} and \Hn{} are the same:
$$dim(Z)=\sum_{\lambda \vdash n}\frac{ n_{\lambda}!}{m_{\lambda}} ,$$
where $n_{\lambda}$ is the number of even parts in $\lambda$ and $m_{\lambda}=\Pi_j( i_j!).$

	\subsection{A Basis for the Center}
		Finding the basis for the center of \Hn{} corresponding to our basis for \AT{} is more complex than in the case of the Nilcoxeter algebra.  In that previous case, it was easy to see what linear combination of diagrams would complement any diagram representing a given basis element of \AT.  In other words, given a basis element $\overline{b}\in \AT$, we needed to find an element $z\in Z(\A)$ such that $\epsilon(bz)=1$ for any $b\in \A$ with $\pi(b)=\overline{b}\in \AT$.  This $z$ was simply the sum of the complements of all such $b\in \A$.  

Because of the relation
$$T_i^2=T_i,$$

the coset of an element of \AT{} tends to be much bigger for $\A=\Hn$ than for $\A=\NC$.  To generate all of the elements of this coset, one needs to first generate all diagrams which are equivalent on the \M, just as was done for \NC.  Then, use the above relation to double crossings and again generate all equivalent diagrams.  Repeat this process until all the new diagrams are redundant.  This algorithm generates redundant diagrams which should be ignored.

\begin{figure}[h!]
\begin{center}
\includegraphics[width=7in]{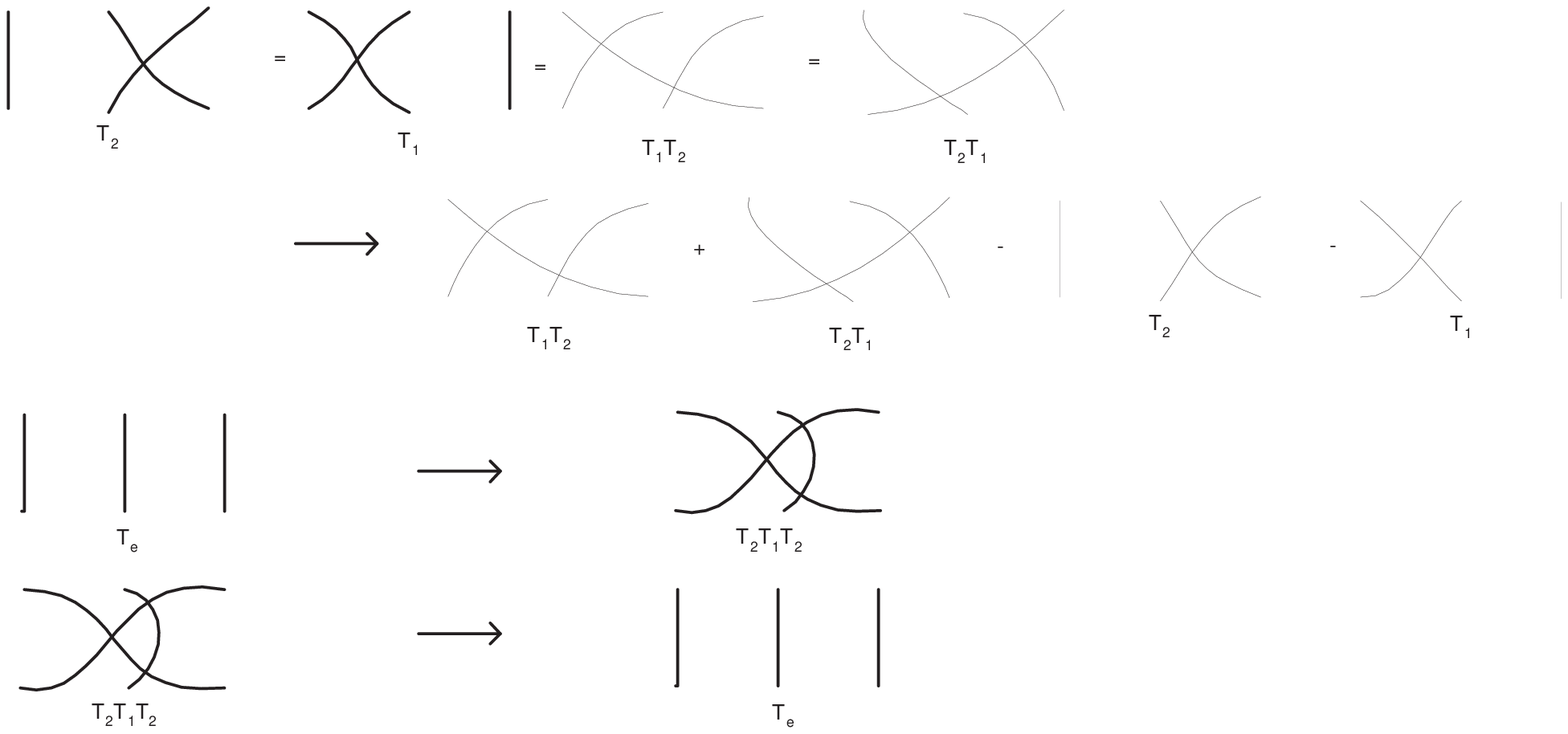}
\end{center}
\caption{This figures shows the cosets corresponding to diagrams on the \M of total thickness three and the basis elements of $Z(\Hn)$ formed by taking the sum of their complements.}
\label{Fig:DualHecke3}
\end{figure}

Then, an element of \Hn{} may have more than one complement. However, note that a diagram in \Hn{} can only have one complement with a given number of crossings.  It then becomes clear that a simple sum of complements like the one we had for \NC{} will not do; one needs a more complicated linear combination of them.  Figure \ref{Fig:DualHecke3} gives the answer for $n=3.$

\bibliography{bib_Centers}
\bibliographystyle{plain}

\end{document}

%% file: macros-II.tex
\usepackage{a4wide,latexsym,amsfonts,amssymb,exscale,enumerate}
\usepackage{amsmath}
\usepackage{amsthm}
\usepackage{hyperref,picins}

\usepackage[enableskew]{youngtab}

\usepackage{pstricks}
\usepackage[tiling]{pst-fill}

\psset{linewidth=0.3pt,dimen=middle}
\psset{xunit=.70cm,yunit=0.70cm}
\psset{arrowsize=1pt 5,arrowlength=0.6,arrowinset=0.7}

\usepackage{graphicx}
\usepackage[all]{xy}
\SelectTips{cm}{}

\newcommand{\NC}{\cal{N}\cal{C}_a}

\newcommand{\refequal}[1]{\xy {\ar@{=}^{#1}
(-1,0)*{};(1,0)*{}};
\endxy}

\usepackage{fancyheadings}
\newcommand{\BOX}{\hbox {$\sqcap$ \kern -1em $\sqcup$}}

\renewcommand{\to}{\rightarrow}

\newcommand{\scs}{\scriptstyle}

\theoremstyle{definition}
\newtheorem{thm}{Theorem}[section]

\newtheorem{lem}[thm]{Lemma}

\newtheorem{prop}[thm]{Proposition}

        \newcommand{\be}{\begin{equation}}
        \newcommand{\ee}{\end{equation}}
        \newcommand{\ba}{\begin{eqnarray}}
        \newcommand{\ea}{\end{eqnarray}}
        \newcommand{\ban}{\begin{eqnarray*}}
        \newcommand{\ean}{\end{eqnarray*}}
        \newcommand{\barr}{\begin{array}}
        \newcommand{\earr}{\end{array}}

\numberwithin{equation}{section}

\def\emph#1{{\sl #1\/}}

\let\phi=\varphi
\let\epsilon=\varepsilon

\usepackage{bbm}
\def\C{{\mathbbm C}}

\def\cal#1{\mathcal{#1}}%
\def\1{\mathbbm{1}}%
\def\shuffle{\,\raise 1pt\hbox{$\scriptscriptstyle\cup{\mskip
               -4mu}\cup$}\,}

%% file: diagrams-quotient.tex
\newcommand{\lowrru}[1]{\xybox{%
  (-8,0)*{};
  (8,0)*{};
  (-6,-18)*{};(6,-9)*{} **\crv{(-6,-13) & (6,-15)} ?(1)*\dir{>};
  (6,-9)*{};(6,0)*{}  **\dir{-} ?(.3)*\dir{ }+(2,0)*{\scs {\bf j}};
}}

\newcommand{\lowllu}[1]{\xybox{%
  (-8,0)*{};
  (8,0)*{};
  (6,-18)*{};(-6,-9)*{} **\crv{(6,-13) & (-6,-15)} ?(1)*\dir{>};
  (-6,-9)*{};(-6,0)*{}  **\dir{-} ?(.3)*\dir{ }+(-2,0)*{\scs {\bf j}};
}}

\newcommand{\bbdl}[1]{\xybox{%
  (2,0);(0,-8) **\crv{(2,-2)&(0,-6)}; ?(.5)*\dir{>}
}}
\newcommand{\bbdlu}[1]{\xybox{%
  (2,0);(0,-8) **\crv{(2,-2)&(0,-6)}; ?(.5)*\dir{<}
}}
\newcommand{\bbdr}[1]{\xybox{%
  (-2,0);(0,-8) **\crv{(-2,-2)&(0,-6)}; ?(.5)*\dir{>}
}}
\newcommand{\bbdru}[1]{\xybox{%
  (-2,0);(0,-8) **\crv{(-2,-2)&(0,-6)}; ?(.5)*\dir{<}
}}

